\documentclass[a4paper,reqno,9pt]{amsart}
\usepackage{amssymb}
\usepackage{cite}
\usepackage{tikz}

\usepackage[left=3cm,right=3cm,top=3cm,bottom=2cm]{geometry} 
\usepackage{amsmath} 

\usepackage{amsfonts}
\usepackage{bigints}
\usepackage{amsmath}
\usepackage{amssymb}
\usepackage{systeme}
\usepackage{mathrsfs}
\usepackage[colorlinks=true]{hyperref}
\numberwithin{equation}{section}
\newtheorem{theorem}{Theorem}[section]
\newtheorem{prop}[theorem]{Proposition}
\newtheorem{lm}[theorem]{Lemma}

\newtheorem{rem}{Remarks}
\newcommand{\R}{\mathbb{R}}

\newcommand{\N}{\mathbb{N}^\ast}

\def\ds{\displaystyle}

\newenvironment{preuve}{{\noindent {\bf Proof. }}}{\hfill {\rule{2.5mm}{2.5mm}}}

\makeatletter
\let\@msm@th@eqref\eqref
\renewcommand{\eqref}[1]{%
	\begingroup
	\leavevmode
	\color{red}%
	\hypersetup{linkbordercolor=[named]{red}}%
	\@msm@th@eqref{#1}%
	\endgroup
}
\makeatother

\author[M.~Amara]{Mustapha Amara}
\address{Department of Mathematics, Faculty of Science of Gab\`es, Research Laboratory Mathematics and Applications LR17ES11; Tunisia}
\email{\sl Mostafa.Amara@fsg.u-gabes.tn}

\title[Global regularity of semi-critical case of $(AQG)_{\alpha,\beta}$ equations in Sobolev spaces]
{Global regularity of semi-critical case of anisotropic quasi-geostrophic  equations in Sobolev spaces}

\begin{document}
	\begin{abstract}
	In this paper, we consider the following anisotropic quasi-geostrophic  equations
	\begin{equation}\tag*{$(AQG)_{\alpha,\beta}$}
			\partial_t\theta+ u_\theta.\nabla\theta +\mu|\partial_1|^{2\alpha}\theta+\nu |\partial_2|^{2\beta}\theta=0,\quad
		u_\theta=\mathcal{R}^{\perp}\theta,		
	\end{equation}
	where $\min\{\alpha,\beta\}=\frac{1}{2}$ et $\max\{\alpha,\beta\}\in \left(\frac{1}{2},1\right)$. This equation is a particular case of the equation introduced by Ye (2019) in \cite{YZ}. In
	this paper, we prove that for any initial data $\theta^0$ in the Sobolev
	space  $H^{s}(\R^2)$, $s >1$, the equation $(AQG)_{\alpha,\beta}$ has a global  solution $\theta$ in $C_b(\R^+,H^s(\R^2)).$
	\end{abstract}
	
	

	\maketitle
	\tableofcontents

	
	\section{\bf Introduction}
	In this paper we deal with the following surface quasi-geostrophic  equation  with fractional horizontal dissipation and fractional vertical thermal diffusion:
	\begin{equation}\tag*{$(AQG)_{\alpha,\beta}$}\label{AQG}
	\begin{cases}
		\partial_t\theta+ u_\theta.\nabla\theta +\mu|\partial_1|^{2\alpha}\theta+\nu |\partial_2|^{2\beta}\theta=0,\\
		u_\theta=\mathcal{R}^{\perp}\theta,\\
		\theta|_{t=0}=\theta^0.
	\end{cases}
\end{equation}
Here, $\alpha, \beta \in (0,1)$, and $\mu, \nu > 0$ are fixed real numbers. The unknown $\theta = \theta(t,x)$ is a real-valued function defined on $\mathbb{R}^+ \times \mathbb{R}^2$, with $\theta^0$ as the given initial data. The operator $\mathcal{R}^\perp$ is defined via Riesz transforms by
\begin{align*}
	\mathcal{R}^\perp\theta=\left( -\mathcal{R}_2\theta,\mathcal{R}_1\theta\right)=\left(-\partial_{2}(-\Delta)^{-\frac{1}{2}}\theta,\partial_{1} (-\Delta)^{-\frac{1}{2}}\theta\right),
\end{align*}
The operators $|\partial_1|^{\gamma}$ and $|\partial_2|^{\gamma}$ are given by 
$$\mathcal{F}(|\partial_1|^{\gamma}f)(\xi)=|\xi_1|^{\gamma}\mathcal{F}(f)(\xi),\quad\mathcal{F}(|\partial_2|^{\gamma}f)(\xi)=|\xi_2|^{2\beta}\mathcal{F}(f)(\xi),\quad\forall \xi=(\xi_1,\xi_2)\in \R^2.$$
\\

The system \ref{AQG} is deeply related, in the case when $\gamma=\alpha=\beta$ and $\mu=\nu$,  to the classical dissipation $\ref{QG}$ equation, with its from as follows
		\begin{equation}\label{QG}\tag{QG}
		\begin{cases}
			\partial_t\theta+ u_\theta.\nabla\theta +\mu(-\Delta)^{\gamma}\theta=0,\\
			u_\theta=\mathcal{R}^\perp\theta,\\
			\theta|_{t=0}=\theta^0.
		\end{cases}
	\end{equation}
	
	The equations \ref{AQG} and \eqref{QG} are special cases of the general quasi-geostrophic approximations for atmospheric and oceanic fluid flow with small
	Rossby and Ekman numbers. The first mathematical studies of this equation was carried out in 1994s by Constantin, Majda and Tabak. For more
	details and mathematical and physical explanations of this model we can consult \cite{DC,CP,JP,CP1}.\\
	
	The equation \eqref{QG} has been extensively studied, with significant focus on global well-posedness. In the sub-critical case ($\gamma > \frac{1}{2}$), the theory is well-developed, with global existence and uniqueness for arbitrary initial data established in various function spaces (see \cite{13A,20A,23A,33A,34A}). However, the critical ($\gamma = \frac{1}{2}$) and super-critical ($\gamma < \frac{1}{2}$) cases are more challenging. In the super-critical case, global results are available only for small initial data (see  \cite{3333,5,23,24}). For the critical case, Constantin, Córdoba, and Wu demonstrated global existence in the Sobolev space $H^1$ under a smallness assumption on the $L^\infty$ norm of $\theta^0$ (see \cite{7}). Additional relevant results can be found in \cite{9,14,15,17}. In 2007, Kiselev, Nazarov, and Volberg proved global well-posedness for arbitrary periodic smooth initial data using an elegant modulus of continuity argument (see \cite{16}). \\

However, for equation \ref{AQG}, there has been little attention to the case where $\alpha \neq \beta$. Recently, Ye, Xu, and Wu in \cite{40} proved the global regularity of system \ref{AQG} with $\mu > 0$, $\nu = 0$, $\alpha = 1$ or $\mu = 0$, $\nu > 0$, $\beta = 1$. Subsequently, Ye in \cite{YZ} established the global regularity in $H^s(\R^2)$, $s\geq 2$, when $\alpha,\beta\in (0,1)$ satisfy
	\begin{equation}\label{ZYCondition}
		\beta>\begin{cases}
			\frac{1}{2\alpha+1},&0<\alpha\leq \frac{1}{2}\\ \\
			\frac{1-\alpha}{2\alpha}&\frac{1}{2}<\alpha<1.
		\end{cases}
	\end{equation}	
Additionally, we refer to our result in \cite{MJ}, where we established global regularity in $H^s(\R^2)$, $s > 2-2\min\{\alpha,\beta\}$, for $\alpha, \beta \in (1/2,1)$ (see also \cite{abc}).\\

It is observed that all works on the anisotropic case are based on the situation where $\alpha > \frac{1}{2}$ or $\beta > \frac{1}{2}$. In this context, the order of horizontal dissipation or vertical thermal diffusion implies a stronger regularization effect, which dominates one of the non-linear terms.
On the other hand, when when $\min\{\alpha, \beta\} \leq \frac{1}{2}$, the regularization is insufficient. In this case, the regularization effect is not strong enough to counterbalance the non-linear term. Therefore, we will call this scenario the semi-critical anisotropic case.\\

The primary aim of this paper is to establish global regularity for the case where $\min\{\alpha,\beta\} = \frac{1}{2}$ and $\max\{\alpha,\beta\} > \frac{1}{2}$ for any initial data $\theta^0$ in the inhomogeneous Sobolev space $H^s(\R^2)$, with $s > 1$. Specifically, the main result of this paper is the following global regularity theorem.
	\begin{theorem}\label{theorem2}
		Let $\alpha,\beta$ two real satisfy 
		$$\min\{\alpha,\beta\}=\frac{1}{2},\quad \frac{1}{2}<\max\{\alpha,\beta\}<1.$$
		Let $\theta^0\in H^s(\R^2)$, $s>1$. Then the system \ref{AQG} admits a unique global solution $\theta$ such that 
		\begin{equation}
			\theta\in C_b(\R^+,H^s(\R^2)),\  |\partial_{1}|^\alpha\theta,|\partial_{2}|^\beta\theta\in L^2(\R^+,H^s(\R^2)).
		\end{equation}
	\end{theorem}
	
We outline the main ideas in the proof of this theorem. Since the existence and uniqueness of local smooth solutions can be established through a standard procedure (see the appendix for details). Therefore, it is sufficient to establish a priori bounds for $\theta$ on the interval $[0,T]$ for any fixed $T > 0$.
\begin{rem}
	\begin{enumerate}
		\item[]
		\item The case where $\min\{\alpha,\beta\} = \frac{1}{2}$ and $\max\{\alpha,\beta\} > \frac{1}{2}$ corresponds to a particular case under the conditions established by Ye in \eqref{ZYCondition}.
	\item As $\alpha$ and $\beta$ play a symmetric role, we assume that $\alpha = \min\{\alpha,\beta\}$ and the proof of the theorem \ref{theorem2} is given in Sect. 3 for the case where $\alpha = \frac{1}{2}$ and $\beta > \frac{1}{2}$. The case where $\beta = \min\{\alpha,\beta\}$ follows a similar pattern, so the details are omitted.
	\item For the sake of simplicity, we will set $\mu=\nu=1$ throughout the paper
	\end{enumerate}
\end{rem}
\section{\bf Notations and Preliminary Results}
	In this preparatory section, we shall introduce some functional spaces and prove some elementary	lemmas that will be used in the proof of Theorem \ref{theorem2}
	\subsection{Notations}
	\begin{enumerate}
		\item $C$ denoted all constants that is a generic constant depending only on the quantities specified in the context.
		\item For $p \in  [1,+\infty]$, we denote by $L^p$	the Lebesgue space endowed with the usual norm $\|.\|_{L^p}$.
		\item $\left(\cdot,\cdot\right)_{L^2}$ denotes the usual inner product in the Hilbert space $L^2(\R^2).$
		\item The Fourier transformation in $\R^2$		
		\begin{equation*}
			\mathcal{F}(f)(\xi)=\widehat{f}(\xi)=\int_{\R^2}e^{-ix.\xi}f(x)dx,\quad \xi\in \R^2.
		\end{equation*}
		The inverse Fourier formula is
		\begin{equation*}
			\mathcal{F}^{-1}(f)(x)=(2\pi)^{-2}\int_{\R^2}e^{i\xi.x}f(\xi)d\xi,\quad x\in \R^2.
		\end{equation*}
\item Let $s\in \R$.
	\begin{enumerate}
		\item Inhomogeneous Sobolev Space $H^s(\R^2)$: 
		\begin{align*}
			&H^{s}(\R^2):=\left\{f\in \mathcal{S}'(\R^2); (1+|\xi|^2)^{s/2}\widehat{f}\in L^2(\R^2)\right\},
		\end{align*} with the norm 
		\begin{align*}
			&	\|f\|_{H^s}=\left(\int_{\R^2}(1+|\xi|^2)^s|\widehat{f}(\xi)|^2d\xi\right)^{\frac{1}{2}},
		\end{align*}
		and the scalar product
		\begin{align*}
			&	\left(f,g\right)_{H^s}=\int_{\R^2}(1+|\xi|^2)^s\widehat{f}(\xi)\overline{\widehat{g}(\xi)}d\xi.
		\end{align*}
		\item Homogeneous Sobolev space $\dot{H}^s(\R^2)$: 
		\begin{align*}
			\dot{H}^{s}(\R^2):=\left\{f\in \mathcal{S}'(\R^2);\widehat{f}\in L^1(\R^2)\mbox{ and }|\xi|^s\widehat{f}\in L^2(\R^2)\right\},
		\end{align*} with the norm 
		\begin{align*}
			\|f\|_{\dot{H}^s}=\left(\int_{\R^2}|\xi|^{2s}|\widehat{f}(\xi)|^2d\xi\right)^{\frac{1}{2}},
		\end{align*}
		and the scalar product 
		\begin{align*}
			&	\left(f,g\right)_{\dot{H}^s}=\int_{\R^2}|\xi|^{2s}\widehat{f}(\xi)\overline{\widehat{g}(\xi)}d\xi.
		\end{align*}
	\end{enumerate}
	\item Let $T > 0$, $r \in [1,+\infty]$ and $X$ be a Banach space. We frequently denote the mixed space $L^r([0, T ], X)$ by $L^r_T(X)$, $C([0,T],X)$ by $C_T(X)$.
	\item Let $X$ be a Banach space. The space of bounded continuous functions from $\R^+$ into $X$, denoted  $C_b(\R^+,X)$, is defined as follows
	$$C_b(\R^+,X):=C(\R^+,X)\bigcap L^\infty(\R^+,X).$$
\end{enumerate}
\subsection{Preliminary results}
We recall some fundamental lemmas concerning Sobolev spaces, beginning with the following classical product rules in homogeneous Sobolev spaces.
\begin{lm}[see \cite{BH}]\label{Lemma1}
	Let $s_1$, $s_2$ be two real numbers such that $s_1<1$ and $s_1+s_2>0$. Then, there exists a positive constant $C=C(s_1,s_2)$ such that for all $f,g\in \dot{H}^{s_1}(\R^2)\bigcap \dot{H}^{s_2}(\R^2)$; 
	\begin{equation}
		\|fg\|_{\dot{H}^{s_1+s_2-1}}\leq C(s_1,s_2) \left(\|f\|_{\dot{H}^{s_1}}\|g\|_{\dot{H}^{s_2}}+\|f\|_{\dot{H}^{s_2}}\|g\|_{\dot{H}^{s_1}}\right).
	\end{equation}
	Moreover, if $s_2<1$, there exists a positive constant $C'=C'(s_1,s_2)$ such that for all $f\in \dot{H}^{s_1}(\R^2)$ and $g\in \dot{H}^{s_2}(\R^2)$; 
	\begin{equation}
		\|fg\|_{\dot{H}^{s_1+s_2-1}}\leq C'(s_1,s_2)\|f\|_{\dot{H}^{s_1}}\|g\|_{\dot{H}^{s_2}}.
	\end{equation}
\end{lm}

The following lemma represents the injection between the homogeneous Sobolev spaces and Lebesgue spaces .
\begin{lm}[see \cite{BH}]\label{Lemma2.2}
	Let $p \in [2, +\infty)$ and $\sigma\in [0,1)$ such that
	$$\frac{1}{p}+\frac{\sigma}{2}=\frac{1}{2}.$$
	Then, there is a constant $C > 0$ such	 that
	$$\|f\|_{L^p(\R^d)}\leq C\||\nabla|^{\sigma}f\|_{L^2(\R^d)}.$$
\end{lm}

We also need the following anisotropic interpolation inequalities.
\begin{lm}[see \cite{MA}]\label{Lemma}
	For $s,s_1,s_2\in\R$ and $z\in[0,1]$, the following anisotropic interpolation inequalities hold true for $i=1,2$:
	\begin{align}
		\label{ing1}	&\||\partial_{i}|^{zs_1+(1-z)s_2}f\|_{H^{s}}\leq \||\partial_{i}|^{s_1}f\|_{H^{s}}^z\||\partial_{i}|^{s_2}f\|_{H^{s}}^{1-z},\\
		\label{ing2}	&\||\partial_{i}|^{zs_1+(1-z)s_2}f\|_{\dot{H}^{s}}\leq \||\partial_{i}|^{s_1}f\|_{\dot{H}^{s_1}}^z\||\partial_{i}|^{s_2}f\|_{\dot{H}^{s_2}}^{1-z}.
	\end{align}
\end{lm}
\par We recall the following important commutator and product estimates:
\begin{lm}\label{Lemma4}
	For $s>1$, if $f,g\in \mathcal{S}(\R^2)$ then for any $\alpha\in(0,1)$
	\begin{align}
		\||\nabla|^s(f g)-f|\nabla|^s g\|_{L^2}\leq s2^sC\left(\||\nabla|^{s+\alpha} f\|_{L^{2}}\||\nabla|^{1-\alpha}g\|_{L^{2}}+\||\nabla|^{s-1+\alpha} g\|_{L^{2}}\||\nabla|^{2-\alpha}f\|_{L^{2}}\right).
	\end{align}
\end{lm}
\begin{preuve} 
We start with:
	\begin{align*}
		\||\nabla|^s(f g)-f|\nabla|^s g\|_{L^2}^2&=\int_{\R^2}\left|\mathcal{F}(|\nabla|^s(f g)-f|\nabla|^s g)(\xi)\right|^2 d\xi\\
		&\leq \int_{\R^2}\left(\int_{\R^2}\left||\xi|^s-|\eta|^s\right||\widehat{f}|(\xi-\eta)|\widehat{g}|(\eta)d\eta\right)^2 d\xi.
	\end{align*}
We use the elementary inequality:
	$$\left||\xi|^s-|\eta|^s\right|\leq s2^{s-1}\left(|\xi-\eta|^s+|\eta|^{s-1}|\xi-\eta|\right)$$
Substituting this inequality into our expression, we get:
	\begin{align*}
		\||\nabla|^s(f g)-f|\nabla|^s g\|_{L^2}^2&\leq s^2 2^{2(s-1)} \int_{\R^2}\left(\int_{\R^2}|\xi-\eta|^s|\widehat{f}(\xi-\eta)||\widehat{g}(\eta)|+|\xi-\eta||\widehat{f}(\xi-\eta)||\eta|^{s-1}|\widehat{g}(\eta)|d\eta\right)^2 d\xi\\
		&\leq s^22^{2s}\|f_1g_1\|^2_{L^2}+s^22^{2s}\|f_2g_2\|^2_{L^2},
	\end{align*}
	where
	\begin{align*}
		&\mathcal{F}(f_1)(\xi)=|\xi|^s|\mathcal{F}(f)(\xi)|&& \mathcal{F}(g_1)(\xi)=|\mathcal{F}(g)(\xi)|\\
		&\mathcal{F}(f_2)(\xi)=|\xi||\mathcal{F}(f)(\xi)|&& \mathcal{F}(g_2)(\xi)=|\xi|^{s-1}|\mathcal{F}(g)(\xi)|\\
	\end{align*}
Applying Hölder's inequality and using Lemma \ref{Lemma2.2}, we obtain:
	\begin{align*}
		\||\nabla|^s(f g)-f|\nabla|^s g\|_{L^2}&\leq s2^s\left( \|f_1\|_{L^{\frac{2}{1-\alpha}}}\|g_1\|_{L^{\frac{2}{\alpha}}}+ \|f_2\|_{L^{\frac{2}{\alpha}}}\|g_1\|_{L^{\frac{2}{1-\alpha}}}\right)\\
		&\leq s2^sC(\alpha) \left( \||\nabla|^{s+\alpha}f\|_{L^{2}}\||\nabla|^{1-\alpha}g\|_{L^{2}}+\||\nabla|^{2-\alpha}f\|_{L^{2}}\||\nabla|^{s-1+\alpha}g\|_{L^{2}}\right).
	\end{align*}
	This completes the proof.
\end{preuve}

The following lemma shows that the Riesz transformation is a continuous operator in Lebesgue space
\begin{lm}[See \cite{JN,SM}]\label{Lemma2}
	For any $p\in (1,+\infty)$, there is a constant $C(p)>0$ such that
	\begin{equation}
		\|\mathcal{R}^\perp\theta\|_{L^p}\leq C(p) \|\theta\|_{L^p}.
	\end{equation}
\end{lm}

We also introduce a result in functional analysis that helps us in the proof in Sect. 4.
\begin{lm}[see \cite{BH}]\label{Lm}
	Let $H$ be Hilbert space and $(u_n)_{n\in\N}$ be a bounded sequence
	of elements in $H$ such that
	\begin{equation*}
		u_n\rightharpoonup u\mbox{ in $H$ }\mbox{ and } \limsup_{n\rightarrow +\infty} \|u_n\|\leq \|u\|,
	\end{equation*}
then $\ds\lim_{n\rightarrow +\infty}\|u_n-u\|=0$.
\end{lm}
\section{\bf Proof of Theorem \ref{theorem2}}
In order to complete the proof of Theorem \ref{theorem2}, it is
sufficient to establish a priori estimates that hold for any fixed $T > 0$. We recall the following propositions, which present the basic bounds:
\begin{prop}
	[see \cite{YZ}] \label{prop1}
	Assume $\theta^0$ satisfies the assumptions stated in Theorem \ref{theorem2} and let $\theta$ be the corresponding solution. Then, for any $t > 0$,
	\begin{equation}
		\|\theta(t)\|_{L^2}^2+2\int_0^t\||\partial_1|^\alpha\theta\|_{L^2}^2d\tau+2\int_0^t\||\partial_2|^\beta\theta\|_{L^2}^2d\tau\leq \|\theta^0\|_{L^2}^2.
	\end{equation}
	Moreover, for any $p\in [2,+\infty]$
	\begin{equation}
		\|\theta(t)\|_{L^p}\leq \|\theta^0\|_{L^p}.
	\end{equation}
\end{prop}
The following proposition was proved by the author in \cite{MA} using estimates established in \cite{YZ}. Briefly, under certain conditions for $\alpha$ and $\beta$, the solution remains uniformly bounded in $H^1$ with respect to time. More specifically
\begin{prop}[see \cite{MA}] \label{prop2}
	Assume $\theta^0$ satisfies the assumptions stated in Theorem \ref{theorem2} and let $\theta$ be the corresponding solution. If $\alpha$ and $\beta$ satisfy \eqref{ZYCondition},
	then, for any $t>0,$
	\begin{equation}
		\|\theta(t)\|_{H^1}^2+\int_0^t\||\partial_1|^\alpha\theta\|_{H^1}^2d\tau+\int_0^t\||\partial_2|^\beta\theta\|_{H^1}^2d\tau\leq C(\theta^0),
	\end{equation}	
	where $C(\theta^0)$ is a constant depending on the initial data $\theta^0$.
\end{prop}

With the uniformly global $H^1$-bound of $\theta$, we are now ready to establish the uniformly global $H^s$-estimate of $\theta$ to complete the proof

\begin{preuve}	
Taking $\sqrt{1+|\nabla|^{2s}}$ to the first equation of \ref{AQG} and testing by $\sqrt{1+|\nabla|^{2s}}\theta$, and using the divergence-free condition, we see that
	\begin{align*}
		\frac{1}{2}\frac{d}{dt}\|\theta(t)\|_{H^s}^2+\||\partial_{1}|^\alpha \theta\|_{H^s}^2+\||\partial_{2}|^\beta|\theta\|_{H^s}^2&=-\left(|\nabla|^s(u_\theta.\nabla\theta)-u_\theta|\nabla|^s\nabla\theta,|\nabla|^s\theta\right)_{L^2}.
	\end{align*}
By the Cauchy-Schwarz inequality, we have 
\begin{equation}
\label{3.4}	\left|\left(|\nabla|^s(u_\theta.\nabla\theta)-u_\theta|\nabla|^s\nabla\theta,|\nabla|^s\theta\right)_{L^2}\right|\leq \||\nabla|^s(u_\theta.\nabla\theta)-u_\theta|\nabla|^s\nabla\theta\|_{L^{2}}\||\nabla|^s\theta\|_{L^2}.
\end{equation}
According to Lemma \ref{Lemma4} we have 
	\begin{align*}
		\||\nabla|^s(u_\theta.\nabla\theta)-u_\theta|\nabla|^{s}\nabla\theta\|_{L^2}&\leq C\||\nabla|^{s+\frac{1}{2}}\theta\|_{L^{2}}\||\nabla|^{\frac{3}{2}}\theta \|_{L^{2}}.
	\end{align*}
In the light of the interpolation inequality, we obtain the fact that $alpha=frac{1}{2}<beta$. 
	\begin{align}
	\nonumber	\||\nabla|^{s+\frac{1}{2}}\theta\|_{L^{2}}&\leq C\left(\||\partial_1|^{s+\frac{1}{2}}\theta\|_{L^{2}}+\||\partial_2|^{s+\frac{1}{2}}\theta\|_{L^{2}}\right)\\
	\label{3.5}	&\leq C \left(\||\partial_{1}|^{\alpha}\theta\|_{H^{s}}+\||\partial_2|^{\beta}\theta\|_{H^{s}}\right).	\end{align}
	Moreover
	\begin{align}
\nonumber	\||\nabla|^{\frac{3}{2}}\theta\|_{L^2}&\leq C\left(\||\partial_{1}|^{\frac{3}{2}}\theta\|_{L^2}+\||\partial_{2}|^{\frac{3}{2}}\theta\|_{L^2}\right)\\
	\label{3.6}	&\leq C \left(\||\partial_{1}|^{\alpha}\theta\|_{H^{1}}+\||\partial_2|^{\beta}\theta\|_{H^{1}}\right)
	\end{align}
	The obtained estimates in \eqref{3.4}, \eqref{3.5} and \eqref{3.6} yield
	\begin{align*}
		\left|\left(|\nabla|^s(u_\theta.\nabla\theta)-u_\theta|\nabla|^s\nabla\theta,|\nabla|^s\theta\right)_{L^2}\right|&\leq C\left(\||\partial_{1}|^{\alpha}\theta\|_{H^{s}}+\||\partial_2|^{\beta}\theta\|_{H^{s}}\right)\left(\||\partial_{1}|^{\alpha}\theta\|_{H^{1}}+\||\partial_2|^{\beta}\theta\|_{H^{1}}\right)
		\| \theta\|_{\dot{H}^s}\\
		&\leq \frac{1}{2}\||\partial_{1}|^{\alpha}\theta\|_{H^{s}}^2+\frac{1}{2}\||\partial_2|^{\beta}\theta\|_{H^{s}}^2+C\left(\||\partial_{1}|^{\alpha}\theta\|_{H^{1}}^2+\||\partial_2|^{\beta}\theta\|_{H^{1}}^2\right)
		\| \theta\|_{H^s}^2.
	\end{align*}
	It is then clear that 
	\begin{align*}
		\frac{d}{dt}\|\theta(t)\|_{H^s}^2+\||\partial_{1}|^\alpha \theta\|_{H^s}^2+\||\partial_{2}|^\beta\theta\|_{H^s}^2\leq C\left(\||\partial_{1}|^{\alpha}\theta\|_{H^{1}}^2+\||\partial_2|^{\beta}\theta\|_{H^{1}}^2\right)
		\| \theta\|_{H^s}^2
	\end{align*}
	Gronwall inequality tells us that for any $t\geq 0$
	\begin{align*}
	\|\theta(t)\|_{H^s}^2+\int_{0}^{t}\||\partial_{1}|^\alpha \theta\|_{H^s}^2d\tau+\int_{0}^t\||\partial_{2}|^\beta\theta\|_{H^s}^2d\tau&\leq \|\theta^0\|_{H^s}^2\exp\left[C\int_0^{+\infty}\left(\||\partial_{1}|^{\alpha}\theta\|_{H^{1}}^2+\||\partial_2|^{\beta}\theta\|_{H^{1}}^2\right)d\tau\right]\\
	&\leq\|\theta^0\|_{H^s}^2\exp(C(\theta^0)).
	\end{align*}
So we get that 
$$\theta\in C_b(\R^+,H^s(\R^2)),\ |\partial_{1}|^{\frac{1}{2}}\theta,|\partial_{2}|^{\beta}\theta \in L^2(\R^+,H^s(\R^2)).$$
\end{preuve}
\section{\bf Appendix. Local well-posedness theory of \ref{AQG}}
To ensure thoroughness, this appendix provides the local existence and uniqueness result for equation \ref{AQG} with initial data $\theta^0\in H^s(\R^2)$ for $s>2-2\min\{\alpha,\beta\}$. Specifically, we demonstrate the following local well-posedness result
\begin{prop}\label{Proposition1}
	Let $\alpha,\beta>0$ such that $\min\{\alpha,\beta\}\leq 1/2$ and $s>2-2\min\{\alpha,\beta\}$, and let $\theta^0\in H^s(\R^2)$. Then there exists a positive time $T$, dependent on $\|\theta^0\|_{H^s}$, for which equation \ref{AQG} has a unique solution.
	$$\theta\in C([0,T],H^s(\R^2)),\quad |\partial_{1}|^{\alpha}\theta,|\partial_{2}|^\beta\theta\in L^2([0,T],H^s(\R^2)).$$
\end{prop}
\subsection{\bf Proof of proposition \ref{Proposition1}}
\subsubsection{\bf  Local existence:}
We start by proving the existence of local solution of \ref{AQG} equation in the space 
$$\mathcal{X}_T=\left\{f\in L^\infty_{T}(H^{s}(\R^2))\cap C_{T}(L^2(\R^2))/\ |\partial_{1}|^\alpha f,|\partial_{2}|^\beta f\in L^2_{T}(H^{s}(\R^2))\right\}.$$
$T$ is a positive real number that will be determined later. To do this, we will use Friedrich's method. For this, we define the spectral cut-off as follows:
$$\mathcal{J}_nf=\mathcal{F}^{-1}\left(\xi\mapsto\chi_{B(0,n)}(\xi)\widehat{f}(\xi)\right),$$
where $n\in \N$ and $B(0,n)=\left\{\xi\in \R^2/\ |\xi|\leq n\right\}$ and 
\begin{align*}
	\chi_{B(0,n)}(\xi):=\begin{cases}
		1&\mbox{ if }\xi\in B(0,n)\\
		0&\mbox{ else.}
	\end{cases}
\end{align*}
We consider the following approximate system 	of \ref{B1},		
\begin{equation}\label{B1}\tag*{$(AQG)_n$}
	\begin{cases}
		\partial_t\theta +|\partial_1|^{2\alpha}\mathcal{J}_n\theta +|\partial_2|^{2\beta}\mathcal{J}_n \theta +\mathcal{J}_n\left(\mathcal{J}_n u_\theta .\nabla\mathcal{J}_n  \theta \right)=0,\\
		u_\theta=\mathcal{R}^{\perp}\theta,\\
		\theta_n|_{t=0}=\mathcal{J}_n\theta ^0(x).
	\end{cases}
\end{equation}
We can prove that for any fixed $n$, there exists a unique local solution $\theta _n$ in $C([0, T_n^\ast),H^s(\R^2))$, by applying the Cauchy-Lipschitz theorem. Additionally, we find that $\mathcal{J}_n \theta_n$ is a solution to equation \ref{B1} with the same initial data, due to the property $\mathcal{J}_n^2=\mathcal{J}_n$. Based on the uniqueness of the solution, we conclude that
$\mathcal{J}_n\theta _n=\theta _n$ and $\mathcal{J}_nu_{\theta _n}=u_{\theta _n}.$
Consequently, the approximation system \ref{B1} simplifies to:
\begin{equation}\label{B2}\tag*{$(AQG)'_n$}
		\partial_t\theta_n +|\partial_1|^{2\alpha}\theta_n +|\partial_2|^{2\beta} \theta_n +\mathcal{J}_n\left( u_{\theta_n} .\nabla \theta_n \right)=0
\end{equation}
By the basic energy estimate, we conclude that $\theta_n$ of \ref{B2} satisfies
\begin{align*}
	\|\theta_n(t)\|_{L^2}^2+2\int_{0}^t\||\partial_{1}|^\alpha \theta_n\|_{L^2}^2d\tau+2\int_{0}^t\||\partial_{2}|^\beta \theta_n\|_{L^2}^2d\tau&\leq \|\theta^0\|_{L^2}^2.
\end{align*}
As a result of Picard's extension theorem, the local solution can be extended to a global one. Furthermore, the $H^s$-estimate allows us to derive
\begin{align}
	\nonumber	\frac{1}{2}\frac{d}{dt}\|\theta_n(t)\|^2_{H^s}+\||\partial_1|^\alpha \theta_n\|^2_{H^s}+\||\partial_2|^\beta \theta_n\|^2_{H^s} &\leq \||\nabla|^s(u_{\theta_n}.\nabla\theta_n)-u_{\theta_n}|\nabla|^s\nabla\theta_n\|_{L^2}\|\theta_n\|_{\dot{H}^s}\\
	\nonumber	&\leq C \||\nabla|^{s+\alpha}\theta_n\|_{L^2}\||\nabla|^{1+\alpha}\theta_n\|_{L^2}\||\theta_n\|_{\dot{H}^s}\\
	\nonumber	&\leq C \||\nabla|^{s+\alpha}\theta_n\|_{L^2}^{\frac{s+1+2\alpha}{s+\alpha}}\|\theta_n\|_{L^2}^{\frac{s-1}{s+\alpha}}\||\theta_n\|_{\dot{H}^s}\\
	\nonumber	&\leq C \||\nabla|^{s+\alpha}\theta_n\|_{L^2}^{\frac{s+1+2\alpha}{s+\alpha}}\|\theta_n\|_{\dot{H}^s}\\
	\label{Est1}	&\leq  \frac{1}{2C}\||\nabla|^{s+\alpha}\theta_n\|_{L^2}^2+C\|\theta_n\|_{\dot{H}^s}^{2\frac{s+\alpha}{s-1}}.
\end{align}
Since $\alpha\leq \beta$ and by using the  interpolation inequality, we get
\begin{align}
	\label{Est2}	\||\nabla|^{s+\alpha}\theta_n\|_{L^2}^2\leq C + C \||\partial_{1}|^\alpha\theta_n\|_{H^s}^2+ C \||\partial_{2}|^\beta\theta_n\|_{H^s}^2.
\end{align}
By combining the estimates \eqref{Est1} and \eqref{Est2}, we obtain
\begin{align*}
	\frac{d}{dt}\|\theta_n(t)\|^2_{H^s}+\||\partial_1|^\alpha \theta_n\|^2_{H^s}+\||\partial_2|^\beta \theta_n\|^2_{H^s} 
	&\leq C\left(1+\|\theta_n\|_{\dot{H}^s}^{2\frac{s+\alpha}{s-1}}\right).
\end{align*}
Define  $T_n=\sup\left\{T>0,\ \|\theta_n\|_{L^\infty_T(\dot{H}^s)}\leq 2\|\theta^0\|_{\dot{H}^s}\right\}$.
Since the function $\left(t\mapsto \|\theta_n\|_{\dot{H}^s}\right)$ is continuous on $\R^+$, the existence of $T_n$ is guaranteed. Moreover, for any $t\in [0,T_n)$, we have 
\begin{align*}
	\frac{d}{dt}\|\theta_n(t)\|^2_{H^s}
	&\leq \tilde{C}\left(1+\|\theta^0\|_{\dot{H}^s}^{2\frac{s+\alpha}{s-1}}\right).
\end{align*}
where $\tilde{C}> 0$ is an absolute constant. One observes that for all $n$
\begin{align*}
	\|\theta_n(t)\|^2_{\dot{H}^s}
	&\leq \|\theta^0\|_{\dot{H}^s}^2+ \tilde{C}t\left(1+\|\theta^0\|_{\dot{H}^s}^{2\frac{s+\alpha}{s-1}}\right).
\end{align*}
Denoting
$$T_0=\frac{\|\theta^0\|_{\dot{H}^s}^2}{2\tilde{C}\left(1+\|\theta^0\|_{\dot{H}^s}^{\frac{2s+1}{s-1}}\right)}>0.$$
Therefore $T_0<T_n$, for any $n\in \N$. Moreover we have, for any $t\in [0,T_0]$
\begin{align*}
\sup_{t\in [0,T_0]}	\|\theta_n(t)\|_{\dot{H}^s}
	&\leq 2 \|\theta^0\|_{\dot{H}^s}^2.
\end{align*}
So we can conclude that $(\theta_n)_{n\in \N}$ is uniformly bounded in $C([0,T_0],H^{s}(\R^2))$.\\

Thanks to these bounds, we have the ability to construct the limit $\theta$ of this suite of solution, using classical argument by combining Ascoli’s Theorem and the Cantor Diagonal Process. Therefore, this is enough for us to show that $\theta$ is solution of \ref{AQG}, moreover this solution satisfy
$$\theta\in L^\infty([0,T_0],H^{s}(\R^2))\cap C([0,T_0],H^{s'}(\R^2)),\ s'<s,$$
and
$$|\partial_{1}|^\alpha \theta,|\partial_{2}|^\beta \theta\in L^2([0,T_0],H^{s}(\R^2)).$$

\subsection{Uniqueness:} To this end, we consider two solutions $\theta$
and $\tilde{\theta}$ of \ref{AQG}, emanating
from the same initial data, and belonging to $\mathcal{X}_{T_0}.$ We denote $\omega=\theta-\tilde{\theta}$ and $u_\omega= u_{\theta}-u_{\tilde{\theta}}$, where $u_{\theta}=\mathcal{R}^\perp \theta$ and $u_{\tilde{\theta}}=\mathcal{R}^\perp \tilde{\theta}$. Then we get
\begin{equation}\label{mU}
\begin{cases}
		\partial_t\omega+|\partial_1|^{2\alpha}\omega+|\partial_2|^{2\beta} \omega+u_{\theta}.\nabla\omega+u_\omega.\nabla\tilde{\theta}=0,\\
		\omega|_{t=0}=0.
\end{cases}
\end{equation}
Applying the basic $L^2$-estimate to \eqref{mU} yields
\begin{align*}
	\frac{1}{2}\frac{d}{dt}\|\omega(t)\|_{L^2}^2+\||\partial_{1}|^\alpha\omega\|_{L^2}^2+\||\partial_{2}|^\beta\omega\|_{L^2}^2=-\left(u_{\omega}.\nabla\tilde{\theta},\omega\right)_{L^2} \leq\|u_{\omega}.\nabla\tilde{\theta}\|_{L^2}\|\omega\|_{L^2}.
\end{align*}
By applying the classical product rules in the homogeneous Sobolev space (Lemma \ref{Lemma1}), with $s_1 = \alpha$ and $s_2 = 1 - \alpha$, we obtain
\begin{align*}
\|u_\omega.\nabla\tilde{\theta}\|_{L^2}=\|u_\omega.\nabla\tilde{\theta}\|_{\dot{H}^{s_1+s_2-1}}
	&\leq C \|u_\omega\|_{\dot{H}^{s_1}}\|\nabla\tilde{\theta}\|_{\dot{H}^{s_2}}\\
	&\leq C \||\nabla|^{\alpha} \omega\|_{L^2} \||\nabla|^{2-\alpha}\tilde{\theta}\|_{L^2}.
\end{align*}
The fact that $\alpha\leq\beta$ we have 
\begin{align*}
	\||\nabla|^{\alpha} \omega\|_{L^2}\leq \| \omega\|_{L^2}+\||\partial_1|^\alpha \omega\|_{L^2}+\||\partial_2|^\beta \omega\|_{L^2},
\end{align*}
and 
\begin{align*}
	\||\nabla|^{2-\alpha}\tilde{\theta}\|_{L^2}\leq \||\nabla|^{\alpha}\tilde{\theta}\|_{H^s} \leq  \| \tilde{\theta}\|_{H^s}+\||\partial_1|^\alpha \tilde{\theta}\|_{H^s}+\||\partial_2|^\beta \tilde{\theta}\|_{H^s}.
\end{align*}
Therefore, we obtain
\begin{align*}
	\frac{d}{dt}\|\omega(t)\|_{L^2}^2+\||\partial_{1}|^\alpha\omega\|_{L^2}^2+\||\partial_{2}|^\beta\omega\|_{L^2}^2\leq C \left(\| \tilde{\theta}\|_{H^s}+\||\partial_1|^\alpha \tilde{\theta}\|_{H^s}+\||\partial_2|^\beta \tilde{\theta}\|_{H^s}\right)\left\|\omega\right\|_{L^2}^2.
\end{align*}
The above estimates along with the Gronwall inequality give
$$\omega(t)=0,\quad\forall t\in [0,T_0].$$
This yields the uniqueness of the solution on $\mathcal{X}_{T_0}.$
\subsection{Continuity} Finally, we begin to show the time continuity of
the solution in $H^{s}(\R^2).$\\
First of all we have $\theta\in C([0,T_0],H^{s'}(\R^2))$, $s'<s$, moreover $\theta$ satisfy
\begin{equation}\label{S3}
	\partial_t \theta +|\partial_{1}|^{2\alpha}\theta+|\partial_{2}|^{2\beta}\theta+u_\theta.\nabla\theta=0.
\end{equation}
Let $(s_k)_{k\in \N}$ be a strictly increasing sequence in $(1,s)$ that converges to $s$. For any $k \in \N$, we have
\begin{equation}\label{integre}
	\frac{1}{2}\frac{d}{dt}\|\theta\|_{H^{s_k}}^2+\||\partial_{1}|^\alpha\theta\|_{H^{s_k}}^2+\||\partial_{2}|^\beta\theta\|_{H^{s_k}}^2=-\left(|\nabla|^{s_k}\left(u_\theta.\nabla\theta\right),|\nabla|^{s_k}\theta\right)_{L^2}.
\end{equation}
This implies that for all $0 \leq t_1 < t_2 \leq T_0$
\begin{equation}
	\label{4.28}	\|\theta(t_2)\|_{\dot{H}^{s_k}}^2+2\int_{t_1}^{t_2}\||\partial_{1}|^\alpha\theta\|_{H^{s_k}}^2d\tau+2\int_{t_1}^{t_2}\||\partial_{2}|^\beta\theta\|_{H^{s_k}}^2d\tau =\|\theta(t_1)\|_{\dot{H}^{s_k}}^2-2\int_{t_1}^{t_2}\left(|\nabla|^{s_k}\left(u_\theta.\nabla\theta\right),|\nabla|^{s_k}\theta\right)_{L^2}d\tau.
\end{equation}

For the second term in the right side. The fact that $u_\theta$ has divergence free and $1<s_k<s$, we get 
\begin{align}
	\nonumber	\int_{t_1}^{t_2}\left|\left(|\nabla|^{s_k}\left(u_\theta.\nabla\theta\right),|\nabla|^{s_k}\theta\right)_{L^2}\right| d\tau
	\nonumber	&\leq  \int_{t_1}^{t_2}\||\nabla|^{s_k}\left(u_\theta.\nabla\theta\right)-u_\theta.|\nabla|^{s_k}\nabla\theta\|_{L^2}\||\nabla|^{s_k}\theta\|_{L^2}d\tau\\
	\nonumber	&\leq s_k2^{s_k}C \int_{t_1}^{t_2}\||\nabla|^{s_k+\alpha}\theta\|_{L^2}\||\nabla|^{2-\alpha}\theta\|_{L^2}\|\theta\|_{H^{s}}d\tau\\
	\label{4.29}	&\leq  \tilde{C}\|\theta\|_{L_{T_0}^\infty(H^{s})} \int_{t_1}^{t_2}\left(1+\||\partial_{1}|^\alpha\theta\|_{H^{s}}^2+\||\partial_{2}|^\beta\theta\|_{H^{s}}^2\right)d\tau.
\end{align}
where $\tilde{C}$ is independent of $k$. Therefore, we obtain the following estimates
\begin{equation}
	\|\theta(t_1)\|_{H^{s_k}}^2\leq 	\|\theta(t_2)\|_{H^{s}}^2 +M_0\int_{t_1}^{t_2}\||\partial_{1}|^\alpha\theta\|_{H^{s}}^2d\tau+M_0\int_{t_1}^{t_2}\||\partial_{2}|^\beta\theta\|_{H^{s}}^2d\tau+M_0(t_2-t_1),
\end{equation}
and 
\begin{equation}
	\|\theta(t_2)\|_{H^{s_k}}^2\leq 	\|\theta(t_1)\|_{H^{s}}^2 +M_0\int_{t_1}^{t_2}\||\partial_{1}|^\alpha\theta\|_{H^{s}}^2d\tau+M_0\int_{t_1}^{t_2}\||\partial_{2}|^\beta\theta\|_{H^{s}}^2d\tau+M_0(t_2-t_1),
\end{equation}
where $M_0=C\left(1+\|\theta\|^2_{L_{T_0}^\infty(H^{s})}\right)$.\\

By applying Lebesgue's dominated convergence theorem and the monotone convergence theorem as $k\rightarrow +\infty$, we get
\begin{equation}\label{t2tot1}
	\|\theta(t_1)\|_{H^{s}}\leq 	\|\theta(t_2)\|_{H^{s}} +M_0\int_{t_1}^{t_2}\||\partial_{1}|^\alpha\theta\|_{H^{s}}^2d\tau+M_0\int_{t_1}^{t_2}\||\partial_{2}|^\beta\theta\|_{H^{s}}^2d\tau+M_0(t_2-t_1),
\end{equation}
and 
\begin{equation}\label{t1tot2}
	\|\theta(t_2)\|_{H^{s}}\leq 	\|\theta(t_1)\|_{H^{s}} +M_0\int_{t_1}^{t_2}\||\partial_{1}|^\alpha\theta\|_{H^{s}}^2d\tau+M_0\int_{t_1}^{t_2}\||\partial_{2}|^\beta\theta\|_{H^{s}}^2d\tau+M_0(t_2-t_1).
\end{equation}
Therefore
\begin{equation}
	\ds\limsup_{t_1\rightarrow t_2}	\|\theta(t_1)\|_{H^{s}}\leq 	\|\theta(t_2)\|_{H^{s}}\mbox{ and }\ds\limsup_{t_2\rightarrow t_1}	\|\theta(t_2)\|_{H^{s}}\leq 	\|\theta(t_1)\|_{H^{s}}
\end{equation}
According to Lemma \ref{Lm}, and the fact that $\theta(t_1) \rightharpoonup \theta(t_2)$ as $t_1 \rightarrow t_2$ and $\theta(t_2) \rightharpoonup \theta(t_1)$ as $t_2 \rightarrow t_1$ in $H^s$, we have
 $$\lim\limits_{t_1\rightarrow t_2}\|\theta(t_1)-\theta(t_2)\|_{H^{s}}=0\mbox{ and }\lim\limits_{t_2\rightarrow t_1}\|\theta(t_1)-\theta(t_2)\|_{H^{s}}=0,$$ which implies that $\theta$ is continuous at $t$ in $H^s(\mathbb{R}^2)$ for any $t \in [0, T_0]$. Therefore, we can conclude that $$\theta\in C([0,T_0],H^{s}(\R^2)).$$
\section{\bf Conclusion}
In \cite{YZ}, Ye addressed the global regularity of the \ref{AQG} equation under certain conditions for $\alpha$ and $\beta$ (condition \eqref{ZYCondition}) with initial data in the Sobolev space $H^s(\mathbb{R}^2)$, where $s \geq 2$. He employed global $H^1$-estimates and $H^2$-estimates in his analysis. Our goal is to improve upon this work under the same conditions for $\alpha$ and $\beta$ but considering initial data in the Sobolev space $H^s(\mathbb{R}^2)$ with $s \geq 2-2\min\{\alpha,\beta\}$. Here, $H^{2-2\min\{\alpha,\beta\}}(\mathbb{R}^2) = H^{2 - 2\alpha}(\mathbb{R}^2) \cap H^{2 - 2\beta}(\mathbb{R}^2)$ represents the critical space for the \eqref{QG} equation in the case $\alpha = \beta$ (see \cite{Mi}). This improvement is based on the uniform $H^1$-boundedness of the local solution, which was proven in \cite{MA}.\\

In \cite{MJ}, the case where $\alpha, \beta > \frac{1}{2}$ is discussed in $H^s(\R^2)$, $s>2-2\min\{\alpha,\beta\}$. Moreover, in this paper, we address the case where $\min\{\alpha, \beta\} = \frac{1}{2}$ and $\max\{\alpha, \beta\} > \frac{1}{2}$. We conjecture that global regularity can be achieved for any initial data in $H^s(\mathbb{R}^2)$, with $s \geq 2-2\min\{\alpha,\beta\}$, given the same conditions on $\alpha$ and $\beta$ as in \eqref{ZYCondition}.

\end{document}